\documentclass[12pt]{article}
\usepackage{graphics}
\usepackage{amsmath}
\usepackage{amssymb}
\usepackage{graphicx}
\usepackage{makeidx}
\usepackage{mathrsfs}
\usepackage{amsfonts}
\usepackage[mathscr]{eucal}
\usepackage{amsmath}

\usepackage{graphics}

\DeclareGraphicsExtensions{.eps,.pdf,.jpg,.png}

\def\RR{\vbox {\hbox to 8.9pt {I\hskip-2.1pt R\hfil}}}

\def\pni{\par\noindent}
\def\vsh{\smallskip}

\def\vsp{\vsh\pni} %% ie. \smallskip + \par

\numberwithin{equation}{section}
\begin{document}
%%%%%%%%%%%%%%%%%%%%%%%%%%%%%%%%%%%%%%%%%%%%%%%%%%%%%%%%%%%%%%%%%%%%%%%%
\vskip 0.50truecm
\font\title=cmbx12 scaled\magstep2
\font\bfs=cmbx12 scaled\magstep1
\font\little=cmr10
\begin{center}
{\title A note on the Lambert W function:\\
Bernstein and Stieltjes properties \\
for a creep model in Linear Viscoelasticity}
%{\title On the creep function of Lambert type \\ in linear viscoelasticty}
%\\[0.25truecm]
%{\title the Jeffreys--Lomnitz law of creep}
 \\  [0.25truecm]
 Francesco MAINARDI$^{(1)}$
Enrico MASINA$^{(2)}$,
  and
 Juan Luis GONZ{\'A}LEZ-SANTANDER$^{(3)}$
%Departamento de Matemáticas. Universidad de Oviedo
%Daniele RITELLI$^{(3)}$
%
\\  [0.25truecm]
$\null^{(1)}$ {\little Department of Physics, University of Bologna, and INFN}
\\ {\little Via Irnerio 46, I-40126 Bologna, Italy}
\\{\little   E-mail: francesco.mainardi@bo.infn.it}
\\[0.25truecm]
$\null^{(2)}$ {\little Department of Mathematics, University of Bologna,}
\\ {\little Piazza Porta San Donato, I-40127 Bologna, Italy}
\\{\little    E-mail: enrico.masina3@unibo.it}
\\[0.25truecm]
$\null^{(3)}$
{\little Departamento de Matem{\'a}ticas, Universidad de Oviedo}
\\{\little C/ Leopoldo Calvo Sotelo 18, 33007 Oviedo, Spain}
\\ {\little E-mail:  gonzalezmarjuan@uniovi.es}
%\null^{(3)}$ {\little Department of  Statistical Sciences,
% University of Bologna},
% \\ {\little Via Belle Arti, I-40126 Bologna, Italy}
% \\{\little E-mail: daniele.ritelli@unibo.it}
%}
\\[0.25truecm]
{\bf Version of
 \today}
 \\[0.25truecm]
 {\bf Paper published in Symmetry (MDPI) 2023, 15, 1654 (13 pages)}
 \\ {\bf DOI: 103390/sym15091654}
\end{center}
%%%%%%%%%%%%
\begin{abstract}
%% Insert your abstract here.
\noindent
%TO CHECK THE ABSTRACT  ACCORDINGLY THE NEW TITLE
The purpose of this note is to propose an application of the Lambert $W$ function in linear viscoelasticity based on the Bernstein and Stieltjes properties of this function. In particular, we recognize the role of its main branch, $W_0(t)$, in a peculiar model of creep with two spectral functions in frequency that completely characterize the creep model. In order to calculate these spectral functions, it turns out that the conjugate symmetry property of the Lambert $W$ function along its branch cut on the negative real axis is essential. We supplement our analysis  by computing  the corresponding relaxation function and providing the plots of all computed functions

\end{abstract}
%
% \maketitle

\vsp {\it 2010 Mathematics Subject Classification (MSC)}:
% 26A33, 33E12,  44A10.
%% 26A33,  %%%%  (main);    Fractional derivatives and integrals
%% 33E12, %% Mittag-Leffler type functions
	33-00,  	%% General reference works (handbooks, dictionaries, bibliographies, % etc.) pertaining to special functions
	34M30, % 	Asymptotics and summation methods for ordinary differential % % equations in the complex domain
42B10,  %%	Fourier and Fourier-Stieltjes transforms and other transforms of ^Fourier type
 44A10,  %% Laplace Transforms
 74Dxx.	%%	Materials of strain-rate type and history type, other materials with %memory (including elastic materials with viscous damping, various viscoelastic % materials)
%\vsp{\it 2010 Physics and Astronomy Classification Scheme (PACS)}:
%% 77.84.s %% Dielectric materials, 77.84.-s
%77.22.Gm. %%% Dielectric relaxation, 77.22.Gm
\vsp

%%%%%%%%%%%%%%%%%%%%%%%%%%%%%%%%%%%%%%%%%%
%\begin{document}

%%%%%%%%%%%%%%%%%%%%%%%%%%%%%%%%%%%%%%%%%%
%\setcounter{section}{-1} %% Remove this when starting to work on the template.
\section{Introduction}
The Lambert function $W(z)$ is defined as the root of the
transcendental {equation} .

\begin{equation}
W(z)  \exp\left(W(z)\right) = z\,.
\label{W-definition}
\end{equation}

The mathematical history of the $W$ function goes back to the 18th century as  outlined in the seminal paper \cite{Corless et al ACM1996} where the interested reader can be informed  about the  analytical and numerical properties of this function.
Further  details can be found in several papers, see, e.g.,
\cite{Jeffrey-et-al MS1996),Kheyfits FCAA2004,Kalugin PhD2011,Kalugin-Jeffrey CR2011,Kalugin-Jeffrey SPRINGER2010},
 Section 4.13 of the handbook \cite{NIST 2010},
 and in the recent book \cite{MezoBOOK2022} with references therein.
For our purposes,  we outline the paper \cite{Kalugin-et-al ITSF2012}.    %%% NEW!!!
 The applications are found in many areas of applied science as outlined,
 e.g., in \cite{Valluri-et-al CJP2000,Jordan CM2014,MezoBOOK2022} and references therein.
  Let us outline also the applications in probability, see, e.g.,
 \cite{Pakes JMAP2011,Vinogradov CSTM2013,Pakes SPL2018}, which could be considered and possibly revised in view of the present
 results in linear viscoelasticity.

In our analysis, we restrict our attention to the main branch $W_0(z)$ of
the $W$ function on the real semiaxis $z=t\ge 0$.
We claim that the Lambert function $W_0(t)$ can be a candidate involved as a possible model for dimensionless creep compliance (i.e., a material function) in linear viscoelasticity. Indeed, according to a previous idea of \cite{GrossBOOK1953}, the property of being a
Bernstein function (that is, with a completely monotonic derivative) for such  material function is proved to be fundamental
as stated in Chap. 2 of the book \cite{MainardiBOOK2022} and references therein. In addition, the property of also being a Stieltjes function
provides further results for this viscoelastic model that, as far as we know, are new in the framework of the literature of linear viscoelasticity.

The plan of this paper is as follows.
In Section \ref{sec2}, we summarize the essentials of linear viscoelasticity in order to point out the properties of the rate of creep of our model based on the Lambert function.
As a consequence, the rate of creep, being a completely monotonic (CM) function with the additional requirement to be a Stieltjes function,
turns out to be expressed in terms of two different spectral functions.
In Section \ref{sec3}, we carry out a numerical and analytical analysis of our new model based on the properties of the Lambert function. We also illustrate with plots the related quantities which can better characterize the model itself. The consistence of our results are validated with MATHEMATICA.
In the conclusions section, we summarize our final remarks.
We have added two appendices: A  to illustrate an alternative method for the calculation of the spectral functions, B to calculate the relaxation function. 

\section{Essentials of Linear Viscoelasticity}\label{sec2}

%We recall that in the linear theory of viscoelasticity, based on the hereditary theory by Volterra, a viscoelastic body is characterized by two distinct but interrelated material functions,
%causal in time (i.e., vanishing for $t<0$):
%the creep compliance $J(t)$ (the strain response to a unit step function of stress), and the relaxation modulus $G(t)$ (the stress response to a unit step function of strain).
%For more details, see e.g.,
%\cite{Christensen BOOK1982},
%\cite{Pipkin BOOK1986},
%\cite{Tschoegl BOOK1989}
% and the recent book \cite{MainardiBOOK2022}.
 
According to the theory of linear viscoelasticity, viscoelastic bodies are characterized by two different interrelated functions, called material functions, which are causal in time (i.e., null for $t<0$). The first material function is the creep compliance $J(t)$, defined as the strain response to a unit step function of stress. The second material function is the relaxation modulus $G(t)$, which is defined as the stress response to a unit step function of strain.
For more details, see, e.g.,
\cite{Christensen BOOK1982,Pipkin BOOK1986,Tschoegl BOOK1989}
 and the recent book \cite{MainardiBOOK2022}.

%By taking $J(0^+)=J_0 >0$ so that $G(0^+)= G_0 =1/J_0$,
%the body is assumed to exhibit a non-vanishing instantaneous response both in creep and relaxation tests. As a consequence, we find it convenient to introduce two dimensionless quantities $\psi(t)$ and $\phi(t)$ as follows:
Taking $J(0^+)=J_0 >0$ (thus $G(0^+)= G_0 =1/J_0$), a viscoelastic body is assumed to show a nonvanishing instantaneous response both in creep and relaxation tests. Therefore, it is convenient to introduce the following dimensionless functions $\psi(t)$ and $\phi(t)$:
\begin{equation}
\label{J-G-definitions}
J(t)= J_0[1+q\, \psi(t)]\,, \quad G(t) = G_0\, \phi(t)\,,
\end{equation}
where $\psi(t)$ is a non-negative increasing function with
\begin{equation}
\psi(0) =0, \label{psi(0)=0},
\end{equation}
and $\phi(t)$ is a non-negative decreasing function with
\begin{equation}
\phi(0)=1. \label{phi(0)=1}
\end{equation}

The nondimensional quantity $q$  takes into account a suitable scaling of the strain, according to convenience, in experimental rheology.
However, in the following we assume for simplicity $q=1$ without loss of generality.

Viscoelastic bodies may be distinguished in solidlike and fluidlike, whether $J(+\infty)$ is finite or infinite, so that
$G(+\infty)= 1/J(+\infty)$ is nonzero or zero, respectively.

As pointed out in most treatises on linear viscoelasticity,
(e.g., in \cite{Pipkin BOOK1986,Tschoegl BOOK1989,MainardiBOOK2022}), the relaxation modulus $G(t)$ can be derived from the corresponding creep compliance $J(t)$  through the Volterra integral equation of the second kind
\begin{equation}
G(t)= \frac{1}{J_0}-\frac{1}{J_0}\,
\int_0^t \!\frac{dJ}{dt'}\, G(t-t')\, dt'.
\label{integral-equation}
\end{equation}

Therefore, according to \eqref{J-G-definitions}, the nondimensional relaxation function $\phi(t)$ obeys the Volterra integral equation
\begin{equation}
\label{integraleqphi}
\phi(t)=1-q\, \int_0^t\frac{d\psi}{dt'}\, \phi(t-t')\, dt'\,, \quad q=1.
\end{equation}

In the {Appendix B} %%% \ref{appendixa},} %MDPI: 
%We all added format like ``Appendix A'' under our rules, please check if some of %them should be ``Appendix B''.
 we calculate the solution of (\ref{integraleqphi}) by using the Laplace transform. Hereafter, we use the following notation for the Laplace transform pair:
\begin{equation}
\label{juxtaposition}
f(t) \div \widetilde f(s) = \mathcal{L}[f(t); s] := \int_0^\infty \! e^{-st}\, f(t)\, dt,
\end{equation}
where the sign $\div$ denotes the juxtaposition of the function $f(t)$ with its Laplace transform $\widetilde{f}(s)$. In order to distinguish the Laplace transform, we use the same notation $f$ for the original function, but overlined with a tilde and with the proper argument $s$.

According to this notation, the solution of (\ref{integraleqphi}) derived in Appendix B %\ref{appendixa} 
is written as%
\begin{equation}
\phi \left( t\right) =\mathcal{L}^{-1}\left[ \frac{1}{s\left( 1+s \,%
\widetilde{\psi} \left( s\right) \right) };t \right] . \label{Phi(t)_general}
\end{equation}

%In linear viscoelasticity, it is quite common to require the existence of positive retardation and relaxation spectra for the material functions $J(t)$ and $G(t)$, as pointed out in the monograph \cite {GrossBOOK1953}  on the mathematical structure of the theories of viscoelasticity. This implies (as formerly proved in \cite{Molinari 1973} and then revisited by \cite{Hanyga 2005}, see also \cite{MainardiBOOK2022}) that $J(t)$ and $G(t)$, and consequently the functions $\psi(t)$ and $\phi(t)$,
%turn out to be Bernstein and CM functions, respectively.

It is quite usual to require the existence of positive retardation and relaxation spectra for the material functions $J(t)$ and $G(t)$, as pointed out in the monograph \cite{GrossBOOK1953}. This implies (as formerly proved in \cite{Molinari 1973} and then revisited by \cite{Hanyga 2005}, as well as in \cite{MainardiBOOK2022}) that $J(t)$ and $G(t)$ and consequently, the dimensionless functions $\psi(t)$ and $\phi(t)$,
turn out to be Bernstein and completely monotonic functions, respectively.

%Here we recall that a completely monotonic (CM) function $f(t)$ is a non-negative, infinitely derivable function with derivatives alternating in sign for $t>0$, like $\exp(-t)$ Also, a Bernstein function is a non-negative function whose derivative is CM, like $1 - \exp(-t)$.
%A necessary and sufficient condition to be a CM function is provided by the Bernstein theorem according to which $f(t)$ is the Laplace transform of a non-negative real function.
%(For more details on  these mathematical properties  the interested reader is referred to the excellent monograph \cite{Schilling-et-al BOOK2012}). According to the Bernstein theorem, we write the rate of creep as follows:

We recall that a completely monotonic (CM) function is a non-negative, infinitely derivable function whose derivatives alternate in sign for $t>0$. Moreover, a Bernstein function is a non-negative function whose derivative is CM.
According to Bernstein's theorem, $f(t)$ is the Laplace transform of a non-negative function, if and only if $f(t)$ is a CM function 
(the interested reader is referred to the excellent monograph \cite{Schilling-et-al BOOK2012} for more details on the mathematical properties of CM functions). Apply Bernstein's theorem to write the rate of creep as follows:
\begin{equation}
\label{CM-spectra}
\psi^\prime(t) = \int_0^\infty \!\!
e^{-rt}\, K(r)\, dr =
\int_0^\infty \!\! e^{-t/\tau}\, H(\tau)\, d\tau \,,
\end{equation}
where  $ K(r)$ and  $H(\tau)$ are  both non-negative functions and denote the required spectra in frequency ($r$) and in time ($\tau=1/r$), respectively.

Then, from (\ref{CM-spectra}), the time spectrum can be determined using the transformation $\tau= 1/r$, so that
\begin{equation}
\label{H(tau)}
H(\tau) =
\frac{K(1/\tau)}{\tau^2}\,.
\end{equation}

We recognize that the Laplace transform of the rate of creep is the iterated Laplace transform of the frequency spectrum, that is, the Stieltjes transform of $K(r)$. Therefore, the Titchmarsh formula provides the inversion of the Stieltjes transform,
see, e.g., \cite{Widder BOOK1946}.
Indeed, since $\psi^\prime (t)$ is CM, taking into account (\ref{CM-spectra}) and (\ref{Laplace_derivative}) with (\ref{psi(0)=0}),
we have for its Laplace transform:
\begin{eqnarray*}
\label{Stieltjes-1p}
\mathcal{L}[\psi^\prime(t);s]&=&
 s\,\widetilde \psi(s)= {\displaystyle \int_0^\infty\! e^{-st}\, \psi^\prime(t)\, dt}\\
&=&{\displaystyle \int_0^\infty\! e^{-st} \left(\int_0^\infty e^{-rt}\, K(r)\, dr
\right) \, dt},
\end{eqnarray*}
and exchanging the order of integration, we  finally obtain:
 \begin{equation}
\label{Stieltjes-1f}
s\,\widetilde \psi(s)=
{\displaystyle \int_0^\infty\left( \int_0^\infty e^{-t(r+s)}\,dt \right)\, K(r)\, dr}=
{\displaystyle \int_0^\infty \frac{K(r)}{r+s}\, dr}.
\end{equation}

We thus recognize that  $K(r)$ is the inverse of the Stieltjes transform of $s\,\widetilde \psi(s)$.
Under suitable conditions, the inversion can be obtained by the Titchmarsh formula. (This formula is found without names in
\cite{Widder BOOK1946}. However, in some papers and books, it is cited with several names, just Titchmarsh in \cite{MainardiBOOK2022},
Stieltjes--Perron in \cite{Henrici BOOK1977},
Gross--Levi in \cite{Apelblat BOOK2011},
Bobylev--Cercignani  in \cite{Aghili-Masomi KJM2014}.
 In any case, it results as a simple exercise in complex analysis as stated in \cite{Gorenflo-Mainardi CISM1997} for the Mittag-Leffler function).
  Thus,
   \begin{equation}
\label{K(r)}
K(r) = \pm \frac{1}{\pi}\,
\left.
\Im[s\,\widetilde\psi (s)]
\right|_{{\displaystyle s =r} e^{\mp i\pi}}\,.
\end{equation}

In addition, if $\psi^\prime(t)$ is a Stieltjes function 
(i.e., the Laplace transform of a CM function), then $K(r)$ turns out to be the Laplace transform of a non-negative function that we denote by $\rho(u)$. 
Indeed,

 \begin{equation}
 \label{Stieltjes-2f}
 \begin{array}{lll}
\psi^\prime(t) = &{\displaystyle \int_0^\infty e^{-r\,t}K(r)\, dr}=
   {\displaystyle \int_0^\infty e^{-rt}
   \left(\int_0^\infty e^{-ru} \rho (u) \, du\right) dr} \\
   & ={\displaystyle \int_0^\infty e^{-r u}
   \left(\int_0^\infty e^{-rt} \, dr \right) \rho(u)\, du}\\
   &= {\displaystyle \int_0^\infty \left(\int_0^\infty
   e^{-r(t+u)}\, d\tau\right) \,\rho(u)\, du}=
   {\displaystyle \int_0^\infty \dfrac{\rho(u)}{t +u} du},
   \end{array}
   \end{equation}
where we have applied
\begin{equation}
\label{K(r)-rho(u)}
K(r) = \int_0^\infty e^{-ru} \rho(u) du.
\end{equation}

Therefore, from (\ref{Stieltjes-2f}), and applying the Titchmarsh formula, we obtain
\begin{equation}
\label{rho(u)-Stieltjes}
\rho(u) = \pm \frac{1}{\pi}\,
\left.
\Im[\psi^\prime(t)]
\right|_{{\displaystyle t =u} e^{\mp i\pi}}\,.
\end{equation}

{Hence, in the case that the rate of creep $\psi^\prime(t)$ is a Stieltjes function, we have another spectral function $\rho(u)$ related to the previous spectral function $K(r)$ by a Laplace transformation.}

\section{Application of the Lambert function in Linear Viscoelasticity}\label{sec3}

\begin{figure} %%[H]

\includegraphics[width=0.45\textwidth]{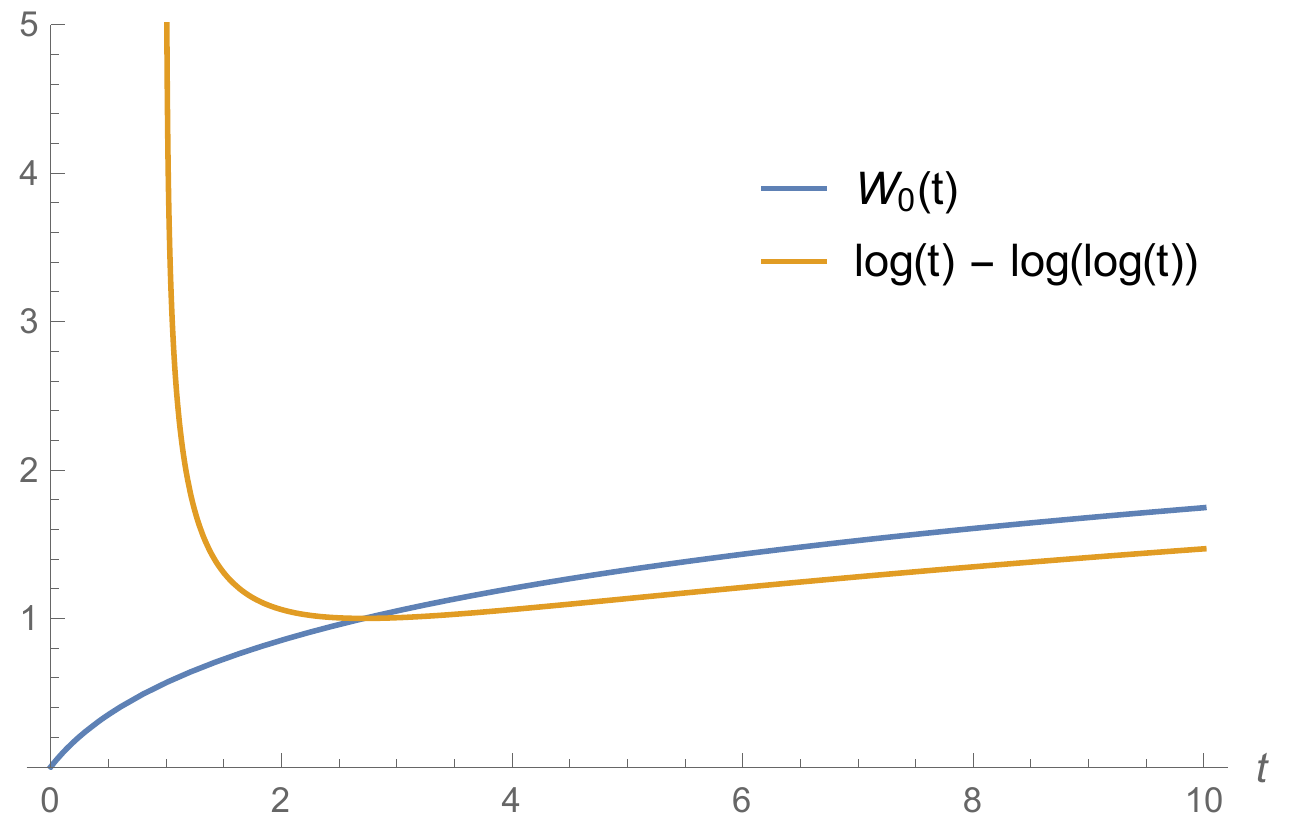}
\includegraphics[width=0.45\textwidth]{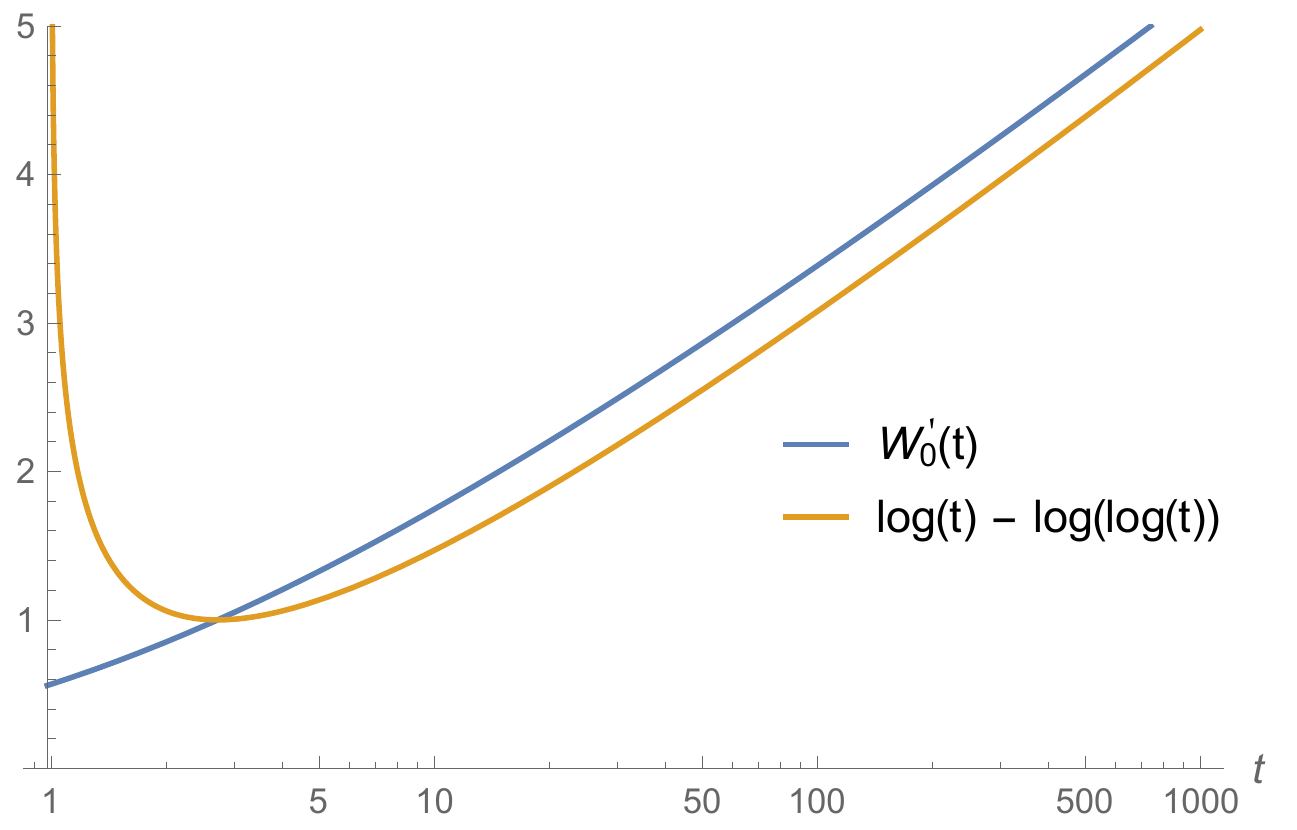}
\caption{ %% Plots of the dimensionless creep function
The  function $W_0(t)$ versus dimensionless time compared with its two-term asymptotic representation.
We have adopted linear scales in the left subfigure for $0<t<10$, and in the right subfigure, log--log scales for $0 < t < 10^3$.
}
\label{fig1}

\end{figure}

After recalling the properties of the Lambert $W$ function, it is worth to see how its plot looks like compared with its asymptotic representation for large arguments. According to Figure \ref{fig1}, the function $W_0(t)$ appears as a positive
one, increasing from 0 at $ t=0$ up to $\infty$ as $t \to \infty$.
The two-term asymptotic representation,
\begin{equation}
W_0(t) \sim \log(t) - \log \log (t), \;t\to \infty,
\end{equation}
fits the function slowly as $t\to\infty$.

We can see that the $W_0(t)$ function is positive increasing up to $\infty$ according to its property of being a Bernstein function,
that is, a non-negative function with a CM first derivative in $t\ge 0$, see, e.g., \cite{Kalugin PhD2011}, p. 116.
Furthermore, its derivative is known to be a Stieltjes function, that is, according to the definition stated in the previous section, a CM function represented by a real Laplace transform of a CM function, where both functions exhibit non-negative spectral functions, see \cite {Kalugin-et-al ITSF2012}.

Differentiating in (\ref{W-definition}) and solving for $W^\prime$, we obtain the following expressions for the derivative of $W_0$:
\begin{equation}
W_0^\prime(t)= \dfrac{1}{e^{W_0(t)} +t}=\dfrac{W_0(t)}{t\left[1+ W_0(t)\right]}\,,\quad t>0\,,
\label{derivative Lambert}
\end{equation}
and its two-term asymptotic representation is
\begin{equation}
W_0^\prime(t) \sim \dfrac{1}{t} - \dfrac{1}{t + \log(t)}, \, \quad
t \to \infty \,.
\label{derivative Lambert asymptotic}
\end{equation}

\begin{figure} %%[H]

\includegraphics[width=0.6\textwidth]{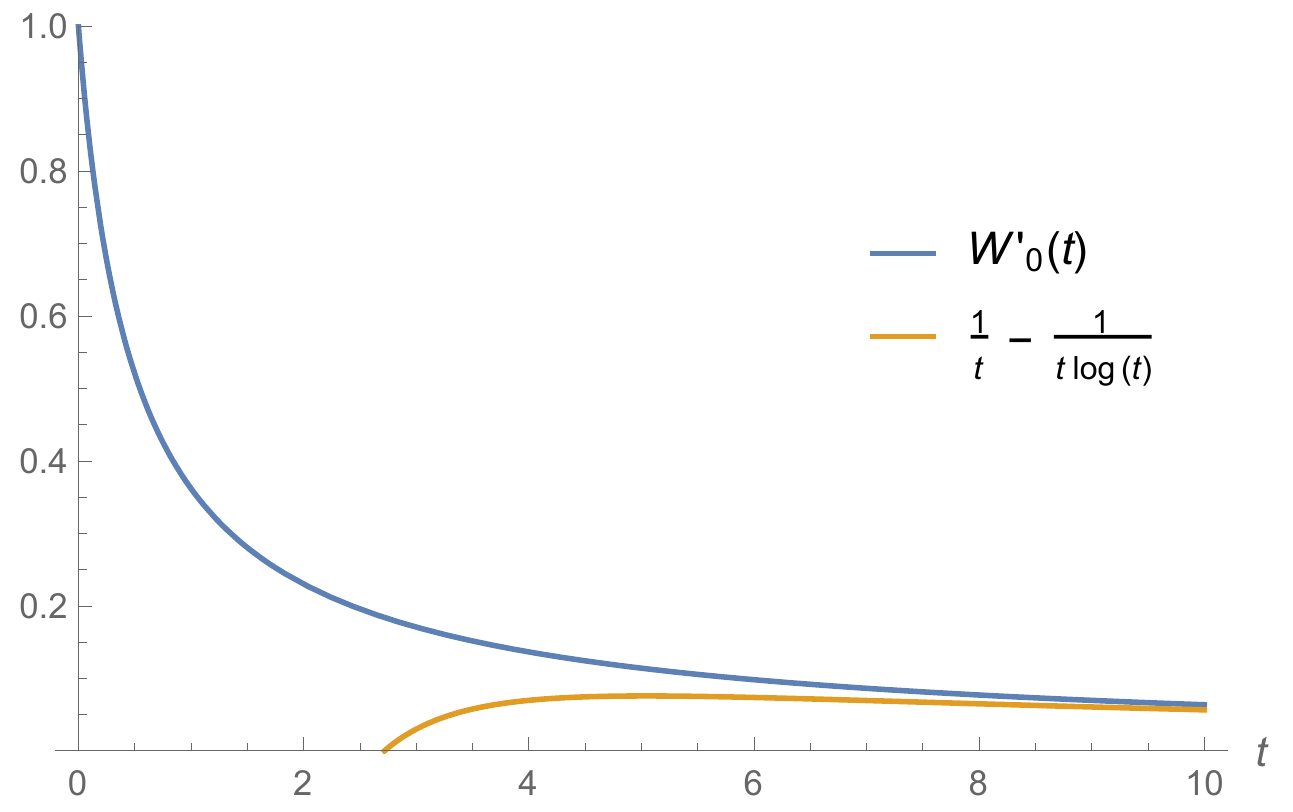}
\caption{The  function $ W_0^\prime(t)$ versus dimensionless time compared with its two-term asymptotic representation for $0<t<10$.}
\label{fig2}

\end{figure}

The plots of the derivative of $W_0(t)$ and its asymptotic representation are depicted in Figure \ref{fig2}.

\vspace{-3pt}

In the proposed viscoelastic model, we assume for the creep that
\begin{equation}
\psi(t) = W_0(t).
\label{psi-Lambert}
\end{equation}

From this assumption, and also taking into account that
 $W_0^\prime\left(t\right)$ is a Stieltjes function, we derive the corresponding spectral functions $K\left(r\right)$ and $\rho\left(u\right)$ presented in the previous section. Indeed, the spectral function $\rho(u)$ can be determined from (\ref{rho(u)-Stieltjes}) and (\ref{psi-Lambert}), i.e.,
\begin{equation}
\rho(u) = -\dfrac{1}{\pi}\,  \Im[W_0^\prime(-u)]\,,\quad u>0,
\label{rho(u)-Lambert}
\end{equation}
where $W_0^\prime(z)$ is computed over or down the negative semiaxis that would be its branch cut in the complex plane.

\begin{figure} %%[H]

\includegraphics[width=0.6\textwidth]{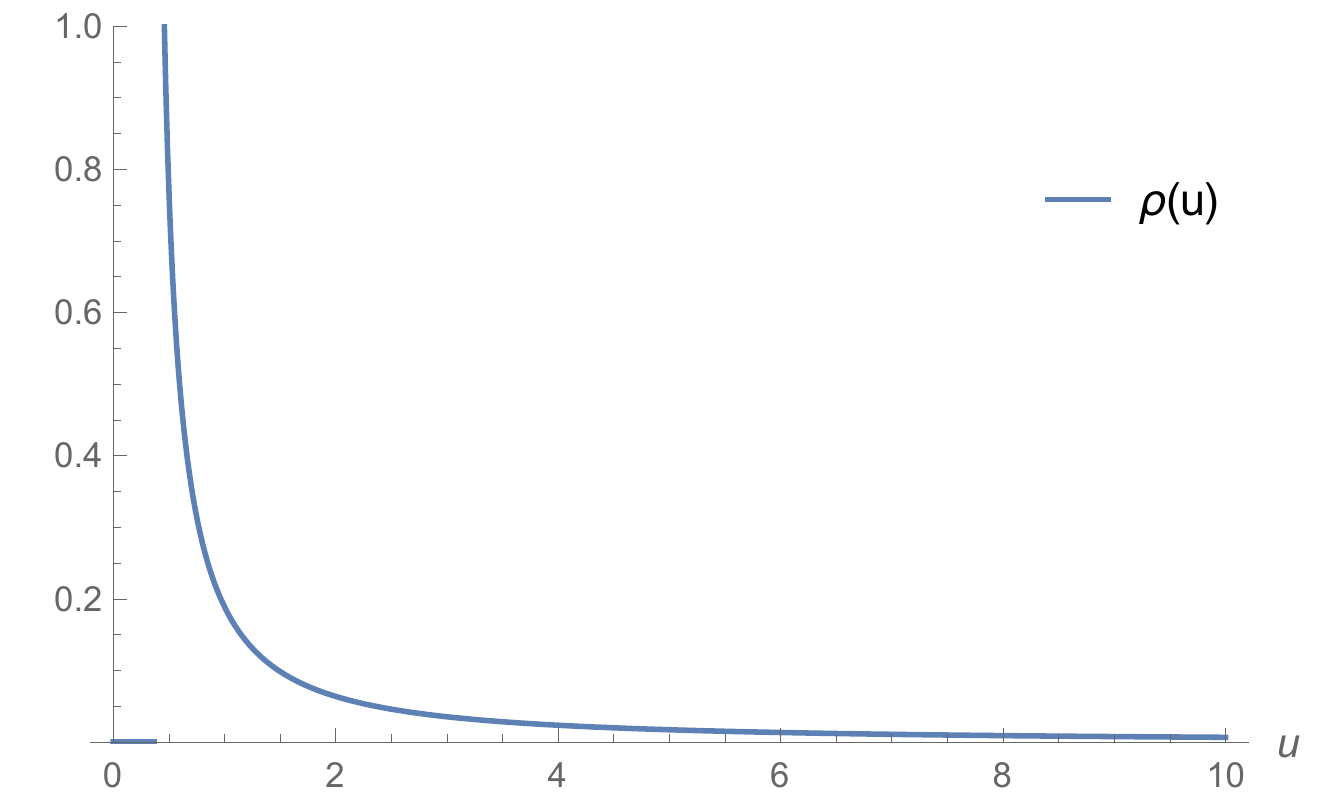}
\caption{The spectral function $\rho(u)$  for $0\le u\le 10 $}
\label{Figure_rho}

\end{figure}

According to the theory of the Titchmarsh formula, the function $\rho(u)$ would be non-negative and continuous.
However, our application of the Titchmarsh formula provides a non-negative function, but discontinuous, see  Figure \ref{Figure_rho}, where we note that $\rho(u) =0$ for $0<u<1/e$ and non-negative for $u>1/e$.
In particular, $\rho(u)$ is decreasing from $\infty$ in the limit $u \to (1/e)^+$ to zero in the limit $u \to +\infty$.

\pagebreak

In order to validate the Stieltjes properties with the related Titchmarsh formula, we need to compute the iterated Laplace transform of the (discontinuous) spectral function $\rho(u)$ and verify that it is the original derivative $W_0^\prime(t)$ (see Equation (\ref{K(r)=L-1[Psi'(t)]})) except for suitable small numerical errors. Indeed, by using MATHEMATICA, this result was verified because
the relative error was small enough, as shown in Figure \ref{fig 4}. Nonetheless, we obtained an alternative derivation of the $\rho\left(u\right)$ function in Appendix A. %\ref{appendixa}. 
It is worth noting that the conjugate symmetry property of the Lambert $W$ function along its branch cut on the negative real axis, i.e., (\ref{symmetry_W}),  turns out to be essential for this derivation.

\begin{figure}

\includegraphics[width=0.8\textwidth]{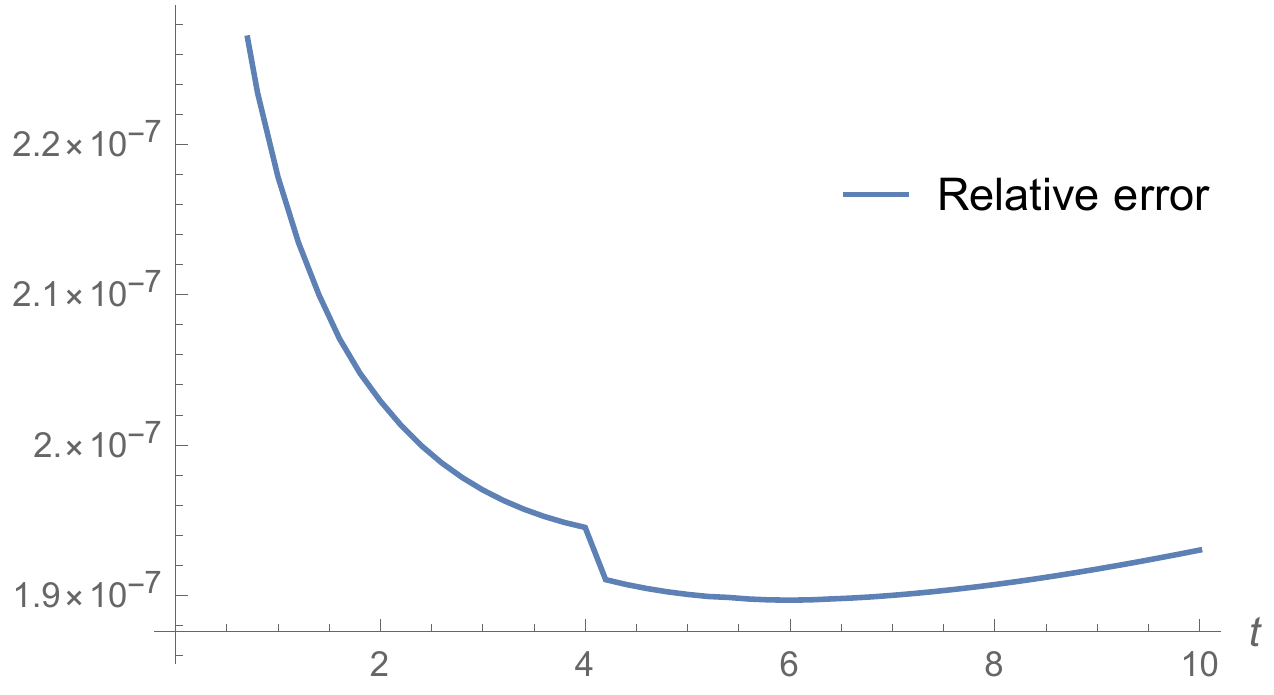}
\caption{Relative error between the iterated Laplace transform of the spectral function $\rho(u)$ and the function $W_0^\prime(t)$
for $t > 0$.}
\label{fig 4}

\end{figure}

We note that a simple reasoning (validated by MATHEMATICA) provides the value of one for the Riemann generalized integral in the positive real semiaxis, that is
\begin{equation}
\int_0^\infty \rho(u) \, du =1.
\label{integral-rho}
\end{equation}

Indeed, taking into account (\ref{rho(u)-Lambert}), performing the change of variables $v =-u$, and knowing that $W_{0}\left(
0\right) =0$, we have%
\begin{eqnarray*}
\int_{0}^{\infty }\rho \left( u\right) du &=&-\frac{1}{\pi }\,\Im%
\int_{0}^{\infty }W_{0}^{\prime }\left( -u\right) du \\
&=&-\frac{1}{\pi } \,\Im\int_{-\infty }^{0}W_{0}^{\prime }
\left( v \right) dv  \\
&=&-\frac{1}{\pi }\lim_{t\rightarrow -\infty }
 \Im\left[ W_{0}\left(
0\right) -W_{0}\left( t\right) \right]  \\
&=&\frac{1}{\pi }\lim_{t\rightarrow -\infty } \Im\left[ W_{0}\left(
t \right) \right] .
\end{eqnarray*}%

According to the property, see \cite{Kalugin-et-al ITSF2012},
$$
\lim_{t\rightarrow -\infty }\Im\left[ W_{0}\left( t\right) \right]
=\pi ,
$$
we conclude that Equation (\ref{integral-rho}) is true.

\begin{figure}

\includegraphics[width=0.8\textwidth]{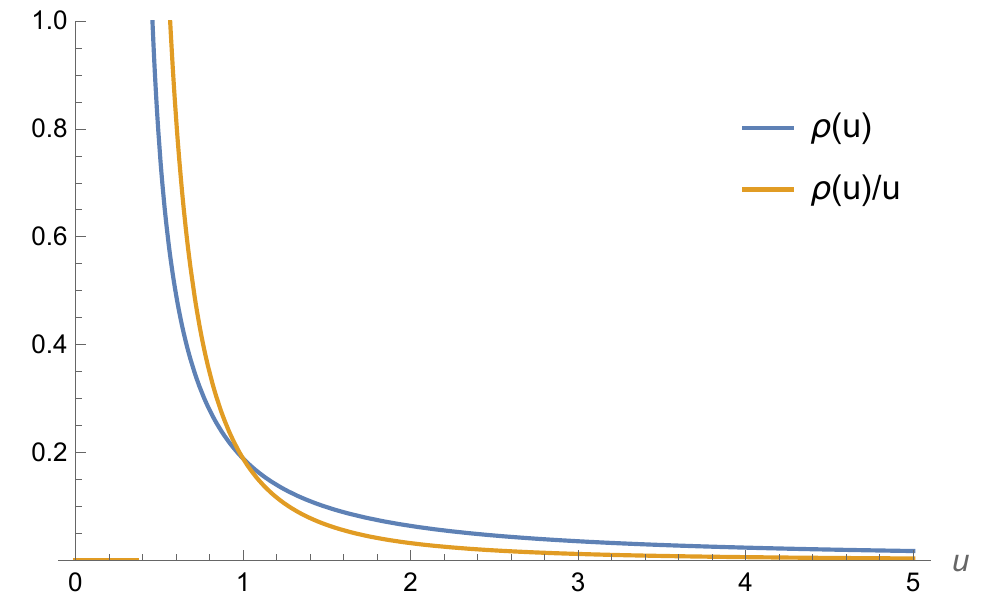}
%% fare new plot!
\caption{
The spectral function $\rho(u)$ compared with the function $\rho(u)/u$ for $0\le u\le 10$.
}
\label{fig4}

\end{figure}

Similar properties are shared by $\rho(u)/u$ for $u\in (0, \infty)$, as validate by MATHEMATICA, but of course with different plots as it is shown in Figure \ref{fig4}.

Indeed, the value of the infinite integral of $\rho(u)/u$ in the positive semiaxis is also one,
\begin{equation}
\int_0^\infty \frac{\rho(u)}{u} \, du =1 = W_0^\prime(0),
\end{equation}
where this result is based on (\ref{Stieltjes-2f}) for $t=0$, (\ref{psi-Lambert}), as well as the fact that $W_0^\prime(0)=1$, according to Figure \ref{fig2}.

\begin{figure}

\includegraphics[width=0.8\textwidth]{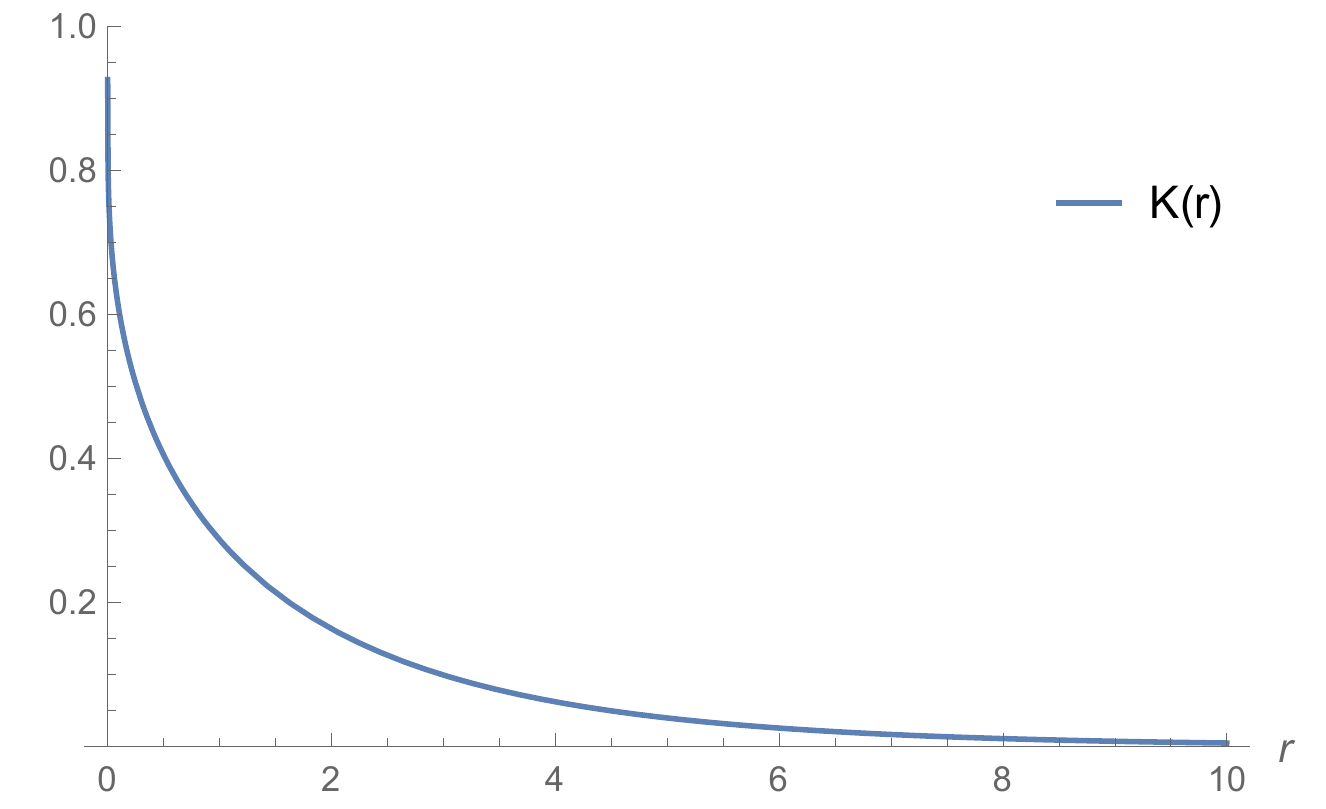}
%% fare new plot!
\caption{
The plot of the spectral function $K(r)$ for $0\le r \le  10$.}
\label{Figure_K(r)}

\end{figure}

According to (\ref{K(r)}) and (\ref{psi-Lambert}), we need the expression of the Laplace transform of $W_0(t)$ in order to calculate the spectral function $K(r)$. However, the latter is not available in the literature. Nevertheless, we can numerically evaluate the spectral function $K(r)$ from (\ref{K(r)-rho(u)}) with the expression given in (\ref{rho(u)-Lambert}) for the $\rho(u)$ function, i.e.,
\begin{equation}
\label{K(r)_Laplace}
K\left( r\right) =-\frac{1}{\pi }\int_{0}^{\infty }\Im \left[ W_{0}^{\prime
}\left( -u\right) \right] du,
\end{equation}
as it is shown in Figure \ref{Figure_K(r)}. Nonetheless, since the $\rho(u)$ function given in (\ref{rho(u)-Lambert}) is subjected to a numerical validation as discussed above, we provide an analytical derivation of (\ref{K(r)_Laplace}) in Appendix A. %%\ref{appendixa}.

We now provide the plot of the spectral function $H(\tau)$ in time.%%
\begin{figure}  %%[H]

\includegraphics[width=0.8\textwidth]{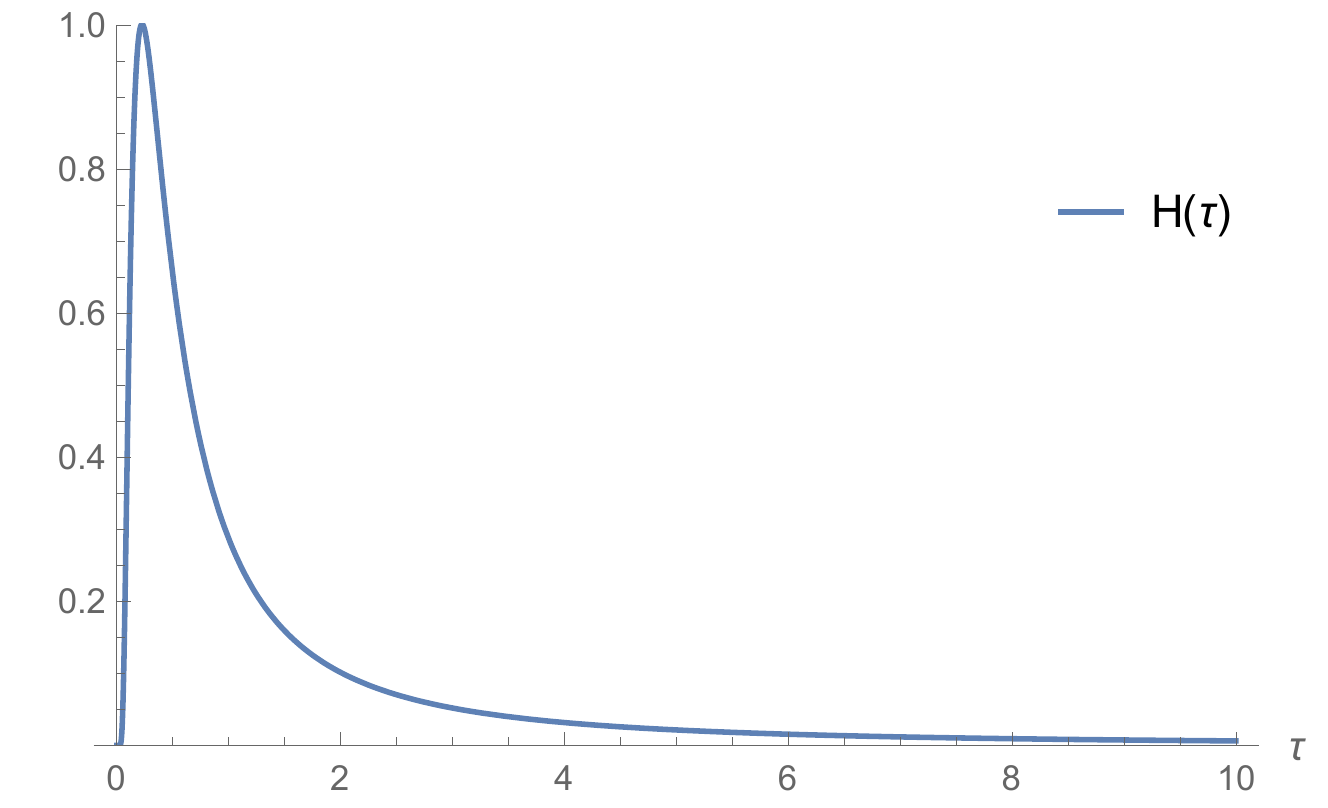}
%% fare new plot!
\caption{
The
 plot of the spectral function $H(\tau)$
for $0\le \tau \le 10$.}

\end{figure}
We recognize that $H(\tau)$ is not a CM function because 
its relation with $K(r)$ in the transformation $r \to \tau= 1/r$ 
cannot preserve this relevant property.

Now, let us calculate the dimensionless relaxation function $\phi(t)$ solving the Volterra integral Equation (\ref{integraleqphi}) for our model, i.e., (\ref{psi-Lambert}). According to (\ref{Phi(t)_general}), we have to compute%
\begin{equation}
\phi \left( t\right) =\mathcal{L}^{-1}\left[ \frac{1}{s\left( 1+s\, \widetilde{W_{0}}\left( s\right)  \right) };t \right] .  \label{Phi_resultado}
\end{equation}

\begin{figure}

\includegraphics[width=0.8\textwidth]{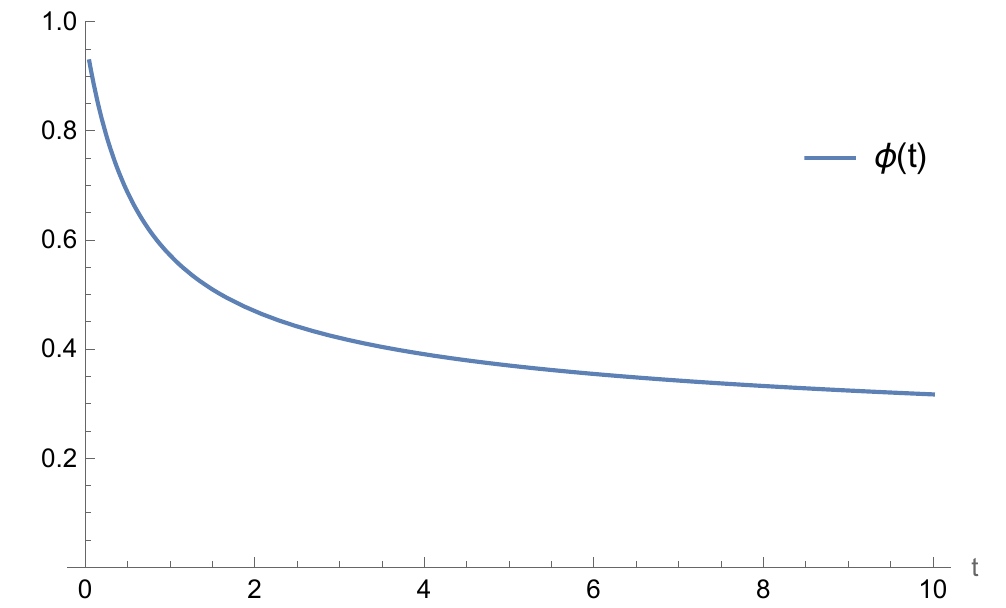} %% fare new plot!

\caption{Dimensionless $\protect\phi \left( t\right) $ function for $t>0$. }
\label{Figure: Phi(t)}
\end{figure}

Figure \ref{Figure: Phi(t)} shows the numerical evaluation of (\ref{Phi_resultado}). Note that the graph of $\phi\left(t\right)$ confirms that it is a non-negative decreasing function with $\phi\left(0\right)=1$, in agreement with (\ref{phi(0)=1}).

%%%%%%%%%%%%%%%%%%%%%%%%%%%%%%%%%%%%%%%%%%
\section{Conclusions}
In this paper, we pointed out some properties of the Lambert function taking advantage of the large existing literature on that function.
In particular, we used these properties to propose a model in linear viscoelasticity that appears to be novel, as far as we know. This model was well characterized, as far as the creep was concerned, with two spectral functions. These spectral functions were numerically validated with MATHEMATICA, as well as analytically derived in Appendix A. %%\ref{appendixa}.

It is worth noting that this model provides an instructive slow-varying creep function, slower than a logarithmic law exhibited by the Lomintz and Becker models known in the geophysical literature, see \cite{MainardiBOOK2022}, Section 2.9.1.

To complete the analysis of our model we provided the calculation of the corresponding relaxation function in Appendix B..

\section*{Acknowledgments }
The research activity of FM
has been carried out in the framework of the activities of the National Group of Mathematical Physics (GNFM, INdAM).
The authors are grateful to PhD R. Garra for useful discussions.

%%%%%%%%%%%%%%%%%%%%%%%%%%%%%%%%%%%%%%%%%%
%\section{Patents}

%This section is not mandatory, but may be added if there are patents resulting from the work reported in this manuscript.

%%%%%%%%%%%%%%%%%%%%%%%%%%%%%%%%%%%%%%%%%%
\vspace{6pt}

%%%%%%%%%%%%%%%%%%%%%%%%%%%%%%%%%%%%%%%%%%
%% optional
%\supplementary{The following supporting information can be downloaded at:  \linksupplementary{s1}, Figure S1: title; Table S1: title; Video S1: title.}


%\begin{tabularx}{\textwidth}{lll}
%\toprule
%\textbf{Name} & \textbf{Type} & \textbf{Description} \\
%\midrule
%S1 & Python script (.py) & Script of python source code used in XX \\
%S2 & Text (.txt) & Script of modelling code used to make Figure X \\
%S3 & Text (.txt) & Raw data from experiment X \\
%S4 & Video (.mp4) & Video demonstrating the hardware in use \\
%\bottomrule
%\end{tabularx}
%}

%%%%%%%%%%%%%%%%%%%%%%%%%%%%%%%%%%%%%%%%%%
\appendix{}

\section{ Alternative method for the calculation of the spectral functions}

\subsection{Inverse Laplace transform of $\protect{W_{0}^{\prime }\left( t\right) }$}

The inverse Laplace transform can be calculated by applying the Bromwich integral:
\begin{equation*}
f\left( t\right) =\mathcal{L}^{-1}\left[ F\left( s\right);t \right]
=\lim_{T\rightarrow \infty }\frac{1}{2\pi i}\int_{\gamma -iT}^{\gamma
+iT}e^{st}F\left( s\right) ds,
\end{equation*}%
where the integration is calculated along the vertical line $\mathrm{Re\,}s=\gamma
$ in the complex plane such that $\gamma $ is greater than the real part of
all singularities of $F\left( s\right) $. In our case, according to (\ref{derivative Lambert}), we have
\begin{equation}
F\left( s\right) =W_{0}^{\prime }\left( s\right) =\frac{W_{0}\left( s\right)
}{s\left[ 1+W_{0}\left( s\right) \right] },  \label{W'0_def}
\end{equation}%
where $W_{0}^{\prime }\left( s\right) $ has two singularities:\ $s_{0}=0$
and $s_{1}=-e^{-1}$ (since $W_{0}\left( -e^{-1}\right) =-1$). Note that $%
s_{0}$ is a simple pole and $s_{1}$ is a branch point. Consider the complex
integral along the contour $C$ depicted in Figure \ref{Figure: integration
path}.
\vspace{-3pt}
\begin{figure} %%[H]
\begin{center}
\includegraphics[width=0.55\textwidth]{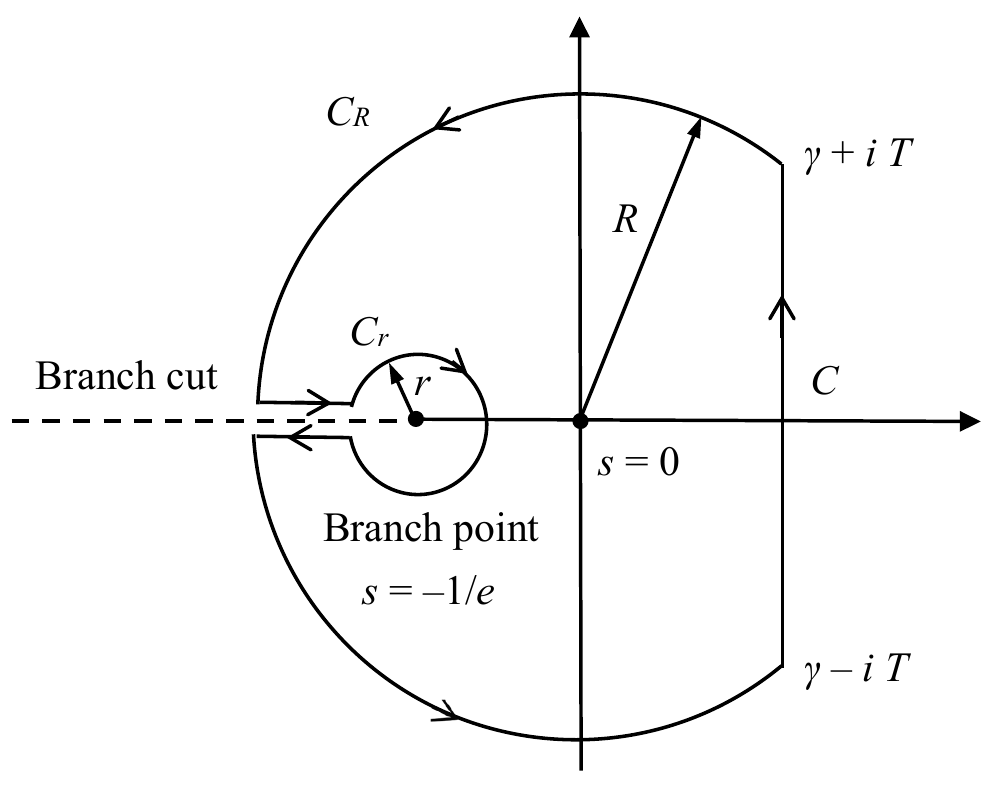}
\caption{Integration path. }
\label{Figure: integration path}
\end{center}
\end{figure}

\begin{eqnarray}
&&\frac{1}{2\pi i}\int_{C}e^{st}W_{0}^{\prime }\left( s\right) ds
\label{int_C}\notag \\
&=&\frac{1}{2\pi i}\left\{ \lim_{T\rightarrow \infty }\int_{\gamma
-iT}^{\gamma +iT}e^{st}W_{0}^{\prime }\left( s\right) ds+\lim_{R\rightarrow
\infty }\int_{C_{R}}e^{st}W_{0}^{\prime }\left( s\right)
ds+\lim_{r\rightarrow 0}\int_{C_{r}}e^{st}W_{0}^{\prime }\left( s\right)
ds\right\}   \\
&&+\frac{1}{2\pi i}\lim_{a\rightarrow 0}\left\{ \int_{-\infty
+a\,i}^{-1/e+a\,i}e^{st}W_{0}^{\prime }\left( s\right)
ds+\int_{-1/e-a\,i}^{-\infty -a\,i}e^{st}W_{0}^{\prime }\left( s\right)
ds\right\} .  \notag
\end{eqnarray}

In order to calculate (\ref{int_C}), let us introduce the following results
\cite{Schiff BOOK1999} {(Chap. 3).} %MDPI: Moved out and changed into standard format.

{\bf Theorem A1}(Cauchy residue theorem)
{\it If $f\left( z\right) $ is analytic within and on a simple, close contour $C$
except at finitely many points $z_{1},z_{2},\ldots ,z_{n}$ lying in the
interior of $C$, then}
\begin{equation*}
\frac{1}{2\pi i}\int_{C}f\left( z\right) dz=\sum_{i=1}^{n}\mathrm{Res}\left(
z_{i}\right) ,
\end{equation*}%
{\it where the integral is taken in the positive direction.}
% \end{Theorem}

{\bf Definition A1}
{\it The function $f\left( z\right) $ has a pole of order $m$ at $z_{0}$, if and
only if}
\begin{equation*}
f\left( z\right) =\frac{h\left( z\right) }{\left( z-z_{0}\right) ^{m}},
\end{equation*}%
{\it where $h\left( z\right) $ is analytic at $z_{0}$ and $h\left( z_{0}\right)
\neq 0$.}
%\end{Definition}

{\bf Theorem A2}(Residues)
{\it If $f\left( z\right) $ has a pole of order $m$ at $z_{0}$, then} %
\begin{equation*}
\mathrm{Res}\left( z_{0}\right) =\frac{1}{\left( m-1\right) !}%
\lim_{z\rightarrow z_{0}}\frac{d^{m-1}}{dz^{m-1}}\left[ \left(
z-z_{0}\right) f\left( z\right) \right] .
\end{equation*}
%\end{Theorem}

Therefore, knowing that $W_{0}\left( 0\right) =0$, we have%
\begin{eqnarray}
\frac{1}{2\pi i}\int_{C}e^{st}W_{0}^{\prime }\left( s\right) ds &=&\mathrm{%
Res}\left( s_{0}\right) =\lim_{s\rightarrow 0}s\,W_{0}^{\prime }\left(
s\right)  \notag \\
&=&\lim_{s\rightarrow 0}\,\frac{W_{0}\left( s\right) }{1+W_{0}\left(
s\right) }=0.  \label{int_C=0}
\end{eqnarray}

Now, consider the following Lemma \cite{Schiff BOOK1999} (Lemma 4.1).

{\bf Lemma A1}
{\it If $F\left( s\right) $ satisfies}
\begin{equation*}
\left\vert F\left( s\right) \right\vert \leq \frac{M}{\left\vert
s\right\vert ^{p}},\quad p>0,
\end{equation*}%
{\it and $C_{R}$ is a circular path of radius $R$ centered at the origin, then,
for $t>0$, we have}
\begin{equation*}
\lim_{R\rightarrow \infty }\int_{C_{R}}e^{st}F\left( s\right) ds=0.
\end{equation*}
%\end{Lemma}

In our case, according to (\ref{W'0_def}), we have
\begin{equation*}
\left\vert F\left( s\right) \right\vert =\frac{\left\vert W_{0}\left(
s\right) \right\vert }{\left\vert s\right\vert \left\vert 1+W_{0}\left(
s\right) \right\vert }\leq \frac{1}{\left\vert s\right\vert },
\end{equation*}%
thus%
\begin{equation}
\lim_{R\rightarrow \infty }\int_{C_{R}}e^{st}W_{0}^{\prime }\left( s\right)
ds=0.  \label{int_CR=0}
\end{equation}

On the other hand, applying the power series of $W_{0}^{\prime }\left(
s\right) $ about the point $s=-e^{-1}$,
\begin{equation}
W_{0}^{\prime }\left( s\right) =\frac{\sqrt{e/2}}{\sqrt{s+1/e}}-\frac{2e}{3}+%
\frac{11e^{3/2}\sqrt{s+1/e}}{12\sqrt{2}}+O\left( s+1/e\right) ,
\label{F(s)_series}
\end{equation}%
and taking $s+\frac{1}{e}=r\,e^{i\theta }$ on $C_{r}$, we have%
\begin{eqnarray*}
&&\lim_{r\rightarrow 0}\int_{C_{r}}e^{st}W_{0}^{\prime }\left( s\right) ds \\
&=&\lim_{r\rightarrow 0}\int_{C_{r}}e^{st}\left[ \frac{\sqrt{e/2}}{\sqrt{%
s+1/e}}-\frac{2e}{3}+\frac{11e^{3/2}\sqrt{s+1/e}}{12\sqrt{2}}+O\left(
s+1/e\right) \right] ds \\
&=&e^{t/e}\lim_{r\rightarrow 0}\int_{-\pi }^{\pi }\exp \left( r\,e^{i\theta
}t\right) \left[ \frac{\sqrt{e/2}}{\sqrt{r}e^{i\theta /2}}-\frac{2e}{3}+%
\frac{11e^{3/2}}{12\sqrt{2}}\sqrt{r}e^{i\theta /2}+O\left( r\,e^{i\theta
}\right) \right] ri\,e^{i\theta }d\theta ,
\end{eqnarray*}%
thus%
\begin{equation}
\lim_{r\rightarrow 0}\int_{C_{r}}e^{st}W_{0}^{\prime }\left( s\right) ds=0.
\label{int_Cr=0}
\end{equation}

Further, we apply the conjugate symmetry property \cite{Kalugin PhD2011} {(Equation (1.25))}
\begin{equation}
W_{0}^{\prime }\left( \bar{s}\right) =\overline{W_{0}^{\prime }\left(
s\right) },\quad s\notin \left( -\infty ,0\right] , \label{symmetry_W}
\end{equation}%
to calculate%
\begin{eqnarray*}
I &=&\frac{1}{2\pi i}\lim_{a\rightarrow 0}\left\{ \int_{-\infty
+a\,i}^{-1/e+a\,i}e^{st}W_{0}^{\prime }\left( s\right)
ds+\int_{-1/e-a\,i}^{-\infty -a\,i}e^{st}W_{0}^{\prime }\left( s\right)
ds\right\} \\
&=&\frac{1}{2\pi i}\lim_{a\rightarrow 0}\left\{ \int_{-\infty
+a\,i}^{-1/e+a\,i}e^{st}\left[ \Re\,W_{0}^{\prime }\left( s\right)
+i\,\Im\,W_{0}^{\prime }\left( s\right) \right] ds\right. \\
&&+\left. \int_{-1/e-a\,i}^{-\infty -a\,i}e^{st}\left[ \Re\,%
W_{0}^{\prime }\left( s\right) -i\,\Im \,W_{0}^{\prime }\left(
s\right) \right] ds\right\} \\
&=&\frac{1}{2\pi i}\left\{ \int_{-\infty }^{-1/e}e^{st}\left[ \Re\,%
W_{0}^{\prime }\left( s\right) +i\,\Im\, W_{0}^{\prime }\left(
s\right) \right] ds\right. \\
&&-\left. \int_{-\infty }^{-1/e}e^{st}\left[ \Re\,W_{0}^{\prime
}\left( s\right) -i\,\Im \, W_{0}^{\prime }\left( s\right) \right]
ds\right\} \\
&=&\frac{1}{\pi }\int_{-\infty }^{-1/e}e^{st}\,\Im \, W_{0}^{\prime
}\left( s\right) ds.
\end{eqnarray*}

Performing the change of variables $s\rightarrow -s$ and taking into account that
$W_{0}^{\prime }\left( s\right) \in
%TCIMACRO{\U{211d} }%
%BeginExpansion
\mathbb{R}
%EndExpansion
$ for $s\geq -1/e$, we obtain
\begin{eqnarray}
I &=&\frac{1}{\pi }\,\Im \, \int_{1/e}^{\infty }e^{-st}W_{0}^{\prime
}\left( -s\right) ds  \label{Int_resultado} \\
&=&\frac{1}{\pi }\,\Im \, \int_{0}^{\infty }e^{-st}W_{0}^{\prime
}\left( -s\right) ds=\frac{1}{\pi }\,\Im \,\mathcal{L}\left[
W_{0}^{\prime }\left( -s\right) \right] .  \notag
\end{eqnarray}

We collect the results given in (\ref{int_C=0}), (\ref{int_CR=0}), (\ref%
{int_Cr=0})\ and (\ref{Int_resultado}) and substitute them into (\ref{int_C}%
)\ to arrive at%
\begin{equation}
\mathcal{L}^{-1}\left[ W_{0}^{\prime }\left( s\right) \right] =\frac{1}{2\pi
i}\lim_{T\rightarrow \infty }\int_{\gamma -iT}^{\gamma
+iT}e^{st}W_{0}^{\prime }\left( s\right) ds=-\frac{1}{\pi }\,\Im%
\mathcal{L}\left[ W_{0}^{\prime }\left( -s\right) \right] .
\label{Inv_Laplace_0}
\end{equation}

Finally, we apply the following result \cite{Schiff BOOK1999} {(Theorem 2.7).}

{\bf Theorem A3.}
{\it Let $f$ be continuous on $\left( 0,\infty \right) $ and of exponential order
$\alpha $. If $f^{\prime }$ is piecewise continuous on $\left[ 0,\infty
\right) $, then}
\begin{equation}
\mathcal{L}\left[ f^{\prime }\left( t\right) \right] =s\,\mathcal{L}\left[
f\left( t\right) \right] -f\left( 0\right) . \label{Laplace_derivative}
\end{equation}
%\end{Theorem}

Therefore, in our case, since $W_{0}\left( 0\right) =0$, we finally arrive
at
\begin{equation}
\mathcal{L}^{-1}\left[ W_{0}^{\prime }\left( s\right);t \right]
 =\frac{t}{\pi }\Im \,\mathcal{L}\left[ W_{0}\left( -s\right);t %
\right]  .  \label{Inv_Laplace}
\end{equation}

\subsection[\appendixname~\thesubsection]{Calculation of $\protect{K(r)}$ and $\protect{\rho(u)}$}

On the one hand, we write (\ref{K(r)-rho(u)}) as
\begin{equation}
\label{Laplace_K(r)}
\psi ^{\prime }\left( t\right) =\mathcal{L}\left[ K\left( r\right) \right].
\end{equation}

On the other hand, knowing that a Stieltjes representation can be expressed as an iterated Laplace transform (\cite{Widder BOOK1946}, p. 325), we rewrite (\ref{Stieltjes-2f}) as

\begin{equation}
\label{Laplace-iterated}
\psi ^{\prime }\left( t\right) =\mathcal{L}\left\{ \mathcal{L}\left[ \rho
\left( u\right) \right] \right\} .
\end{equation}

From (\ref{Laplace_K(r)}) and (\ref{Laplace-iterated}), we obtain

\begin{equation}
\label{K(r)_inverse_Laplace}
K\left( r\right) =\mathcal{L}\left[ \rho \left( u\right) \right] =\mathcal{L}%
^{-1}\left[ \psi ^{\prime }\left( t\right) \right] ,
\end{equation}

which, according to the assumption of our model (\ref{psi-Lambert}), becomes

\begin{equation}
K\left( r\right) =\mathcal{L}\left[ \rho \left( u\right) \right] =\mathcal{L}%
^{-1}\left[ W_{0}^{\prime }\left( t\right) \right] .
\label{K(r)=L-1[Psi'(t)]}
\end{equation}

According to (\ref{Inv_Laplace_0}), we have,%
\begin{equation}
K\left( r\right) =\mathcal{L}\left[ \rho \left( u\right) \right] =-\frac{1}{%
\pi }\,\Im \,\mathcal{L}\left[ W_{0}^{\prime }\left( -u\right) \right]
,  \label{K(r)_resultado}
\end{equation}%
thus
\begin{equation}
\rho \left( u\right) =-\frac{1}{\pi }\,\Im\,W_{0}^{\prime }\left(
-u\right) .  \label{rho(u)_resultado}
\end{equation}

Equations (\ref{K(r)_resultado}) and (\ref{rho(u)_resultado}) are both in agreement with the results given in (\ref{K(r)_Laplace}) and (\ref{rho(u)-Lambert}), respectively, derived in the body of the paper.

\section{Calculation of $\protect\phi \left( t\right) $}
Let us solve (\ref{integraleqphi}), i.e.,%
\begin{equation}
\int_{0}^{t}\psi ^{\prime }\left( \tau \right) \phi \left( t-\tau \right)
d\tau =1-\phi \left( t\right) .  \label{Eq_Volterra}
\end{equation}

In order to calculate (\ref{Eq_Volterra}), let us introduce the following
result \cite{Schiff BOOK1999} {(Theorem 2.39).}

{\bf Theorem B1.}(Convolution theorem)
{\it If $f$ and $g$ are piecewise continuous on $\left[ 0,\infty \right) $ and of
exponential order $\alpha $, then}%
\begin{equation*}
\mathcal{L}\left[ \left( f\ast g\right) \left( t\right) \right] =\mathcal{L}%
\left[ f\left( t\right) \right] \,\mathcal{L}\left[ g\left( t\right) \right]
,\qquad \mathrm{Re}\,s>\alpha ,
\end{equation*}%
{\it where the convolution is given by the integral}%
\begin{equation*}
\left( f\ast g\right) \left( t\right) =\int_{0}^{t}f\left( \tau \right)
g\left( t-\tau \right) d\tau .
\end{equation*}
%\end{Theorem}

Therefore, by applying the Laplace transform to (\ref{Eq_Volterra}), we obtain%
\begin{equation}
\mathcal{L}\left[ \psi ^{\prime }\left( t\right) \right] \,\mathcal{L}\left[
\phi \left( t\right) \right] =\mathcal{L}\left[ 1-\phi \left( t\right) %
\right] .  \label{Laplace_Eq_Volterra}
\end{equation}

We apply (\ref{Laplace_derivative}), taking into account (\ref{psi(0)=0}), i.e.,
\begin{equation*}
\mathcal{L}\left[ \psi ^{\prime }\left( t\right) \right] =s\,\mathcal{L}\left[
\psi \left( t\right) \right] ,
\end{equation*}%
and recall that
\begin{equation*}
\mathcal{L}\left[ 1\right] =\int_{0}^{\infty }e^{-st}dt=\frac{1}{s},
\end{equation*}%
in order to rewrite (\ref{Laplace_Eq_Volterra})\ as%
\begin{equation*}
s\,\mathcal{L}\left[ \psi \left( t\right) \right] \mathcal{L}\left[ \phi
\left( t\right) \right] =\frac{1}{s}-\mathcal{L}\left[ \phi \left( t\right) %
\right] .
\end{equation*}

Solving for $\phi \left( t\right) $ yields%
\begin{equation}
\label{Phi_solution}
\phi \left( t\right) =\mathcal{L}^{-1}\left[ \frac{1}{s\left( 1+s\, \widetilde{\psi} \left( s\right) \right) };t\right] .
\end{equation}

%%%%%%%%%%%%%%%%%%%%%%%%%%%%%%%%%%%%%%%

%\reftitle{References}

%=====================================
% References, variant B: internal bibliography
%=====================================

\end{document}